\documentclass{amsart}
\usepackage{amssymb,latexsym, amsmath, amscd}
\usepackage{enumerate}
\usepackage{multirow}
\usepackage{fullpage, booktabs}
\usepackage{graphicx}
\usepackage[all]{xy}
\makeatletter
\@namedef{subjclassname@2010}{%
\textup{2010} Mathematics Subject Classification}
\makeatother

\newtheorem{theorem}{Theorem}

\theoremstyle{definition}
\newtheorem{definition}[theorem]{Definition}
\newtheorem{remark}[theorem]{Remark}
\newtheorem{example}[theorem]{Example}

\numberwithin{equation}{section}

\def\Z{\mathbb Z}

\def\Y{\Upsilon}
\def\X{\Pi}
\def\R{\mathbb R}

\def\Q{{Q}}
\def\QQn{{Q_n^{\square}}}

\def\PP{\overline{{P}}}

\begin{document}

\title{Round twin groups on few strands}

\author{Jacob Mostovoy}


\maketitle

\begin{abstract} 
We study the space $\Q_n$ of all configurations of $n$ ordered points on the circle such that no three points coincide, and in which one of the points (say, the last one) is fixed.  We compute its fundamental group for $n<6$ and describe its homology for $n=6,7$. For arbitrary $n$, we compute its first homology and its Euler characteristic. 

We use three geometric approaches. On one hand, $\Q_n$ is naturally defined as the complement to an arrangement of codimension-2 subtori in a real torus. On the other hand, $\Q_n$ is homotopy equivalent to an explicit nonpositively curved cubical complex. Finally, $\Q_n$ can also be assembled from no-3-equal manifolds of the real line.

We also observe that, up to homotopy, $\Q_n$ may be identified with a subspace of the oriented double cover of the moduli space $\overline{\mathcal{M}}_{0,n}(\R)$ of stable real rational curves with $n$ marked points. This gives an embedding of $\pi_1\Q_n$ into the pure cactus group. As a corollary, we see that $\pi_1\Q_n$ is residually nilpotent.

\end{abstract}


\section{Introduction}

The spaces $M_{n,k}(\R^m)$ of configurations of $n$ ordered points in $\R^m$, no $k$ of which may coincide, have been thoroughly studied in the literature. While they do not form an operad like the configuration spaces of distinct points in $\R^n$, they do form a bimodule over the little $m$-disk operad, and this allows one to compute their cohomology \cite{DT}. When $m=1$ and $k=3$, they are Eilenberg-MacLane spaces and their topology is defined by the fundamental group. These fundamental groups are known as pure twin groups (also called flat braid groups, or planar braid groups, etc.) and have also been studied in their own right.
  
Not much is known yet about the spaces $M_{n,k}(X)$ of configurations of $n$ ordered points without $k$-fold coincidences on an arbitrary space $X$. In this note, we consider the space $M_{n,3}(S^1)$ of configurations on a circle $S^1$. 

Fix a point $\infty\in S^1$. Then 
$M_{n,3}(S^1)$ is a Cartesian product of $S^1$ with the subspace $$\Q_n\subset M_{n,3}(S^1)$$ of configurations whose last point is $\infty$. Therefore, in what follows, we will concentrate our attention on $\Q_n$.

\medskip

We explore three approaches to the homology of $\Q_n$. First, we consider  $\Q_n$ as the complement of a discriminant in an $n-1$-dimensional torus. This produces a long exact sequence from which we compute the rank of the first homology group of $\Q_n$. Then, we produce an explicit cubical complex with the homotopy type of $\Q_n$. While the (co)homology of this complex is elusive at this point, it does allow one to compute the Euler characteristic of $\Q_n$. Moreover, this cubical complex is readily seen to be CAT(0) and, as a consequence, its topology is determined by the fundamental group. We call this group the \emph{round twin group on $n$ strands} and denote it by $\X_n$. Finally, we relate $\Q_n$ to the spaces $M_{m,3}(\R)$; namely, we show that $\Q_n$ is the classifying space of a certain graph of pure twin groups on $n-1$ and $n-2$ strands. We obtain another long exact sequence which allows us to estimate the cohomological dimension of $\Q_n$.

\medskip

None of the above approaches gives us a complete description of the topology of $\Q_n$ for all $n$, although for $n\leq 7$ we get all the Betti numbers. When  $n\leq 5$, the group $\X_n$ can be identified explicitly. As for $\Q_6$, we exhibit a 2-dimensional cubical complex with the same homotopy type and explore its symmetries.

\medskip

The round twin group $\X_n$ can be thought of as a subgroup of the pure cactus group $\Gamma_n$, the fundamental group of  the moduli space $\overline{\mathcal{M}}_{0,n}(\R)$ of stable real rational curves with $n$ marked points. In fact, $\Q_5$ is homotopy equivalent to the double cover of $\overline{\mathcal{M}}_{0,5}(\R)$ and $\Q_6$  to the complement of a certain link in the double cover of $\overline{\mathcal{M}}_{0,6}(\R)$. This gives an alternative way of describing the topology if $\Q_n$ for $n<7$.

\section{The first Betti number}
Let $$T^{n-1}=\{(x_1, \ldots, x_{n-1})\,|\, x_i\in S^1\}$$
be the $n-1$-dimensional torus and set $x_n=\infty\in S^1$. 
For $1\leq i<j<k\leq n$,  let $A_{ijk}\in T^{n-1}$ be the subset consisting of all the points satisfying
$$x_i=x_j=x_k.$$
Write $$A_n=\bigcup_{1\leq i<j<k\leq n} A_{ijk}.$$
\begin{definition} The space $\Q_n$ is the complement of $A_n$ in $T^{n-1}$. 
\end{definition}
\begin{theorem}\label{firstbetti}
$H_1(\Q_n)$ is a free abelian group of rank $n-1+\binom{n-1}{3}$ for all $n>3$.
\end{theorem}
\begin{proof}
By duality (Theorem~3.44 of \cite{Hatcher}), we have
$$H_i(T^{n-1},\Q_n)=H^{n-i-1}(A_n,\Z),$$
so that the exact sequence of the pair $(T^{n-1}, \Q_n)$ becomes
$$\ldots\longrightarrow   H_{i+1}(T^{n-1})\   \to H^{n-i-2}(A_n,\Z) \to H_i(\Q_n)\to H_i(T^{n-1})\to H^{n-i-1}(A_n,\Z)\longrightarrow\ldots$$
The space $A_n$ is a union of $n-3$-dimensional tori that intersect along tori of smaller dimension. As a consequence, $$H^{n-2}(A_n,\Z)=0$$ and  $H^{n-3}(A_n,\Z)$ is the free abelian group generated by the $\binom{n}{3}$ classes of the $A_{ijk}$.
The map $$H_{2}(T^{n-1})\   \to H^{n-3}(A_n,\Z)$$
in the above long exact sequence is injective with free abelian quotient. 
As a consequence, $H_1(\Q_n)$ is free abelian of rank
$$(n-1)+\binom{n}{3}-\binom{n-1}{2}=n-1+\binom{n-1}{3}.$$
\end{proof}


\section{The structure of a cubical complex on $\Q_n$}\label{cubical}

Call two configurations in $\Q_n$ equivalent if they can be taken into each other by an orientation-preserving homeomorphism 
$S^1\to S^1$ that fixes $\infty$. Equivalence classes of configurations can be encoded by the following combinatorial data that we call a \emph{$k$-fold crossing type}:
\begin{itemize}
\item an integer  $k$ with $0\leq k\leq [n/2]$;
\item a choice of $k$ disjoint pairs of numbers between $1$ and $n$;
\item a cyclic order on the set with $n-k$ elements that consists of the $k$ chosen pairs and those numbers from 1 to $n$ that do not belong to any pair. 
\end{itemize}
This is the same as a choice of $k$ numbers $m_1,\ldots,m_k$ with $1\leq m_1$, $m_j +1 < m_{j+1}$ and $m_k<n$ and a permutation $\sigma$ of the set $\{1,\ldots, n\}$ such that $\sigma(m_j)<\sigma(m_j + 1)$ and $\sigma(n)=n$.
In what follows, we will denote a $k$-fold crossing type by the string of numbers $\sigma(1),\ldots,\sigma(n)$ in which the pairs
$\sigma(m_j)\sigma(m_j + 1)$ are enclosed in square brackets; see example below. 

Identify $S^1-\{\infty\}$ with $\R$ and give $S^1$ a total order by declaring $x<\infty$ for $x\in\R$. Then, the equivalence class of configurations that corresponds to the permutation $\sigma$ and the integers $m_1,\ldots,m_k$ is
$$\{(x_1,\ldots, x_{n-1})\,|\, x_{\sigma(1)}\leq \ldots \leq x_{\sigma(n-1)}\leq x_{\sigma(n)}=x_{n}=\infty,\ \text{and\ } x_{\sigma(i)}=x_{\sigma(i+1)} \Leftrightarrow i=m_j\ \text{for some\ }j    \}.$$
It is an open simplex of dimension $n-1-k$. 

\begin{example}
The equivalence class corresponding to the crossing type $21[34][56]$ consists of the points
$$\{(x_1,\ldots, x_5)\in T^5\,|\, x_2<x_1<x_3=x_4<x_5=\infty\}.$$
\end{example}

A $k$-fold crossing type can be ``resolved'' in  $2^k$ ways to obtain a $(k-1)$-fold crossing type. The resolution consists of removing one of the pairs from the list of pairs and ordering the elements of the eliminated pair in one of the two possible ways\footnote{We should assume here that $n>2$; for $n=2$ both resolutions of the only possible 1-fold crossing type are the same.}. The equivalence class corresponding to a $k$-fold crossing type lies in the boundary of the $2^k$ equivalence classes that correspond to the  $2^k$ resolutions of this crossing type.

The decomposition of  $\Q_n$ into such open simplices is not a CW-complex. Consider the dual cubical complex $\QQn$, which has one cubical $k$-cell for each $k$-fold crossing type, and whose boundary maps are dual to the adjacencies of the equivalence classes in $\Q_n$. It is a CW-complex, and it is homotopy equivalent to $\Q_n$. Indeed, the boundary of each cubical cell of $\QQn$ consists of cells of smaller dimension. On the other hand, the closures of the cells in $\Q_n$ and in $\QQn$ define two covers of the respective spaces with contractible finite intersections and isomorphic nerves so that $\QQn$ is homotopy equivalent to $\Q_n$.

The $k$-cell count in $\QQn$ gives
\begin{theorem}\label{euler}
$$\chi(\Q_n)=(n-1)! \cdot 2^{1-\frac{n}{2}} \cos\left(\frac{\pi n}{4}\right).$$
\end{theorem}
\begin{proof} 
The number of $k$-fold crossing types such that $x_n=n$ is not included in any of the pairs is $\frac{(n-1)!}{2^k}\binom {n-k-1}{k}$; the number of those with $x_n$ included in some pair is $\frac{(n-1)!}{2^{k-1}}\binom {n-k-1}{k-1}$. 
Therefore, the total number of $k$-cubes  in  
$\QQn$ is equal to  $$\frac{n!(n-k-1)! }{2^k k!(n-2k)!}.$$
The alternating sum of these numbers was computed in \cite{GF}, see also \cite{OEISA009014}.
\end{proof}

. 

The universal cover of $\QQn$ is a cubical CAT(0) complex. Gromov's link condition in our case is implied by the following statement: a set of $k$ distinct 1-fold crossing types such that any pair of types appears in the resolution of a 2-fold crossing type actually appears in the resolution of a $k$-fold crossing type. As a consequence, we have:
\begin{theorem}
 $\Q_n$ is an Eilenberg--Maclane space.
\end{theorem}
In particular, the topology of  $\Q_n$ is determined by $\pi_1\Q_n=\X_n$.

\section{$\Q_n$ and the pure twin groups}\label{graphsofgroups}
The space $\Q_n$ naturally contains the space $M_{n-1,3}(\R)$ as the subspace of configurations that have the point at infinity with multiplicity 1. 
\begin{theorem}\label{les}
For $n>2$, there is a long exact sequence
$$ \ldots\leftarrow  \left({H}^k(M_{n-2,3}(\R),\Z)\right)^{n-1}\xleftarrow{d} H^k(M_{n-1,3}(\R),\Z) \xleftarrow{i^*} H^k(\Q_n,\Z)\leftarrow   \left( {H}^{k-1}(M_{n-2,3}(\R),\Z)\right)^{n-1}\leftarrow \ldots,$$
where the map $i^*$ is induced by the inclusion of  $M_{n-1,3}(\R)$ into $\Q_n$. 
\end{theorem}
The most immediate corollary of this is an upper bound on the cohomological dimension of $\Q_n$. 
\begin{theorem}\label{cohdim}
$H^i(\Q_n,\Z)=0$ for all $i>[(n+1)/3]$.
\end{theorem}
Indeed, for  $i>[(n+1)/3]$ we have $H^i(M_{n-2,3}(\R), \Z)=H^i(M_{n-1,3}(\R), \Z)=0$, see \cite{B}.

\begin{proof}[Proof of Theorem~\ref{les}]
Let $n>2$. For brevity, write $M_{n-1}$ for  $M_{n-1,3}(\R)$; the points of the configurations in  $M_{n-1}$ are labelled by the natural numbers from 1 to $n-1$.  
Denote by $M_{n-2}^{j}$ a copy of $M_{n-2,3}(\R)$ whose configurations are  labelled by natural numbers from 1 to $n-1$ with the label $j$ omitted. 

For a configuration $x\in M_{n-2}^{j}$ define 
$\rho_j(x)\in  M_{n-1}$ 
by adding a point $x_j$ with the label $j$ 
to the right of all the points of $x$. Similarly, 
$\lambda_j(x)$ is defined by adding $x_j$ to the left of $x$. 
These concatenation operations can be considered as maps
$$M_{n-2}^{j}\to M_{n-1}.$$
Indeed, in both cases one can assume that all the points of each configuration in $M_{n-2}^{j}$ lie in some fixed open interval and choose $x_j$ to be a fixed point outside of this interval.  

In the union $$M_{n-1}\sqcup \bigsqcup_{1\leq j< n} M_{n-2}^{j} \times [-1,1],$$
identify, for each $x\in M_{n-2}^{j}$, the point $(x,-1)$ with 
$\lambda_j(x)\in M_{n-1}$ and the point $(x,1)$ with $\rho_j(x)\in M_{n-1}$. Denote the resulting space by $\Q_{n}'$.

We now prove that the space $\Q_{n}'$ is homotopy equivalent to $\Q_{n}$. Let $\Q_{n,I}\subset \Q_n$ be the subspace consisting of the configurations which have either one or two points outside of a finite open interval $I\subset \R\subset S^1$, that is, zero or one points in $\R\backslash I$. We claim that the inclusion of $\Q_{n,I}$ into $\Q_n$ is a homotopy equivalence. Indeed, for any pair of finite open intervals $I\subseteq I'$, the inclusion  $\Q_{n,I}\to \Q_{n,I'}$ is a homotopy equivalence. Any compact subspace of $\Q_n$ lies in $\Q_{n,I}$ for some finite interval $I$, and, therefore, the inclusion 
$\Q_{n,I}\to\Q_n$ induces an isomorphism of homotopy groups. Since both $\Q_{n,I}$ and $\Q_n$ have homotopy type of cell complexes, they are homotopy equivalent.

The space $\Q_{n,I}$ is covered by $n$ open subspaces: $U_j$ with $1\leq j < n$ and $V$. The subspace $U_j$  consists of the configurations whose intersection with $S^1- I$ consists of the point $x_j$ and $x_n=\infty$; the subspace $V$ consists of those configurations whose only point at $\infty$ is $x_n$.

$V$ is homotopy equivalent to $M_{n-1}$ while $U_j$ is homeomorphic to $M_{n-2}^j\times [-1,1]$. The sets $U_j$  are disjoint, while 
$U_j\cap V$ can be identified, up to homotopy, with two copies of $M_{n-2}^j$. The inclusion map
$$U_j\cap V \to U_j$$
is equivalent to the inclusion 
$$M_{n-2}^j\times \{-1\}\sqcup M_{n-2}^j\times \{1\} \to M_{n-2}^j\times [-1,1]$$
of the bases into the cylinder. The map
$$U_j\cap V \to V$$
is equivalent to the concatenations
$x\mapsto \rho_j(x)$
on one copy of $M_{n-2}^j$ and to
$x\mapsto \lambda_j(x)$
on the other copy; this shows that $\Q_n'$ is homotopy equivalent to $\Q_n$.

The inclusion map $M_{n-1} \to \Q_n'$ is a cofibration
$$M_{n-1} \to \Q_n'\to \Sigma \left(* \sqcup \bigsqcup_{1\leq j< n} M_{n-2,j}\right),$$
where * is a one-point space and $\Sigma$ denotes the suspension.
This cofibration gives rise to a long exact sequence in cohomology
{\small
$${{\ldots \leftarrow  \widetilde{H}^{k+1}\left(\Sigma \left(* \sqcup \bigsqcup_{1\leq j< n} M_{n-2}^{j}\right),\Z\right)\leftarrow H^k(M_{n-1},\Z) \leftarrow H^k(\Q_n,\Z)\leftarrow   \widetilde{H}^k\left(\Sigma \left(* \sqcup \bigsqcup_{1\leq j< n} M_{n-2}^{j}\right),\Z\right)\leftarrow \ldots,}}$$}
which, in view of the suspension isomorphism, produces the long exact sequence of the theorem.
\end{proof}

\medskip

The preceding construction can be described  in terms of Bass-Serre theory. 
\begin{figure}[ht]
\includegraphics[width=2.5in]{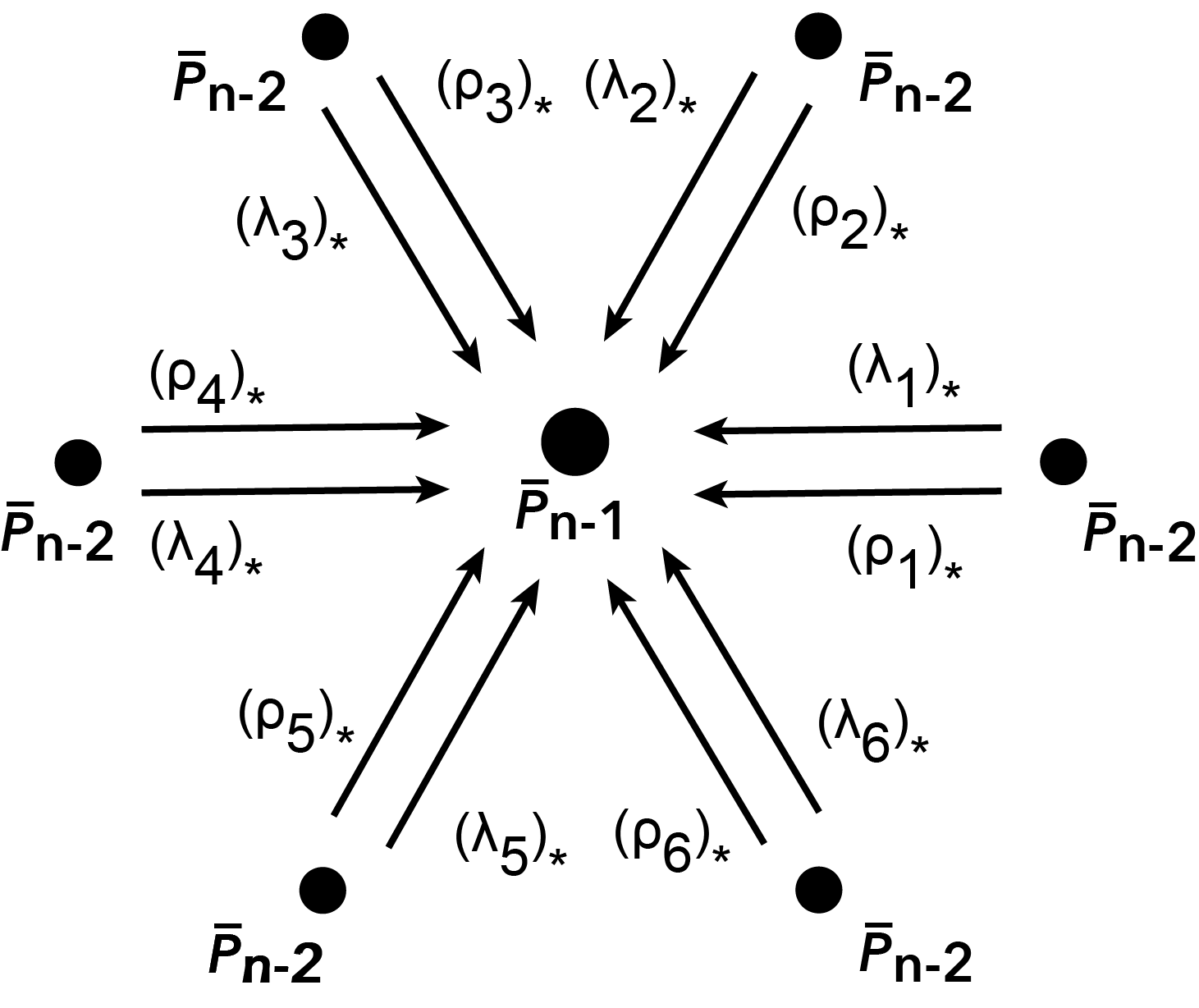}
\caption{The graph of groups $\Phi_n$ for $n=7$.}\label{phisix}
\end{figure}

Choose basepoints $a\in M_{n-2}$ and $b\in M_{n-1}$ and $2n-2$ paths  $g_{j+}$ and $g_{j-}$ which connect $b$ to $\rho_{j}(a)$ and to $\lambda_{j}(a)$ respectively. Then each of the maps $\rho_{j},\lambda_{j}: M_{n-2}\to M_{n-1}$ induces a homomorphism of fundamental groups $\pi_1 M_{n-2}\to \pi_1 M_{n-1}$; denote these homomorphisms by $(\rho_{j})_*$ and 
$(\lambda_{j})_*$ respectively. Recall \cite{Kh2} that the fundamental groups $\pi_1 M_{n-2}$ and $\pi_1 M_{n-1}$ are the \emph{pure twin groups} which we denote by $\PP_{n-2}$ and $\PP_{n-1}$ respectively.

Let $\Phi_n$ be the graph of groups as in Figure~\ref{phisix}. 

\begin{theorem}\label{grgr}
$\Q_n$ has the homotopy type of the classifying space $B\Phi_n$
\end{theorem}
Indeed, the construction of the space $\Q_n'$, which is homotopy equivalent to $\Q_n$ coincides with the construction of the classifying space of this graph, see \cite[Chapter 1.B]{Hatcher}.

\begin{remark}
According to Theorem~1B.11 of \cite{Hatcher}, Theorem~\ref{grgr} also implies that the spaces $\Q_n$ have the homotopy type of the Eilenberg--MacLane spaces $K(\X_n,1)$.
\end{remark}

\section{Round twin groups on few strands}


\begin{theorem}
For $n\leq 4$ the round twin group $\X_n$ is free on $2^{n-1}$ generators.
\end{theorem}
\begin{proof}
The space $\Q_1$ consists of a single point, while $\Q_2$ is a circle and $\Q_3$ is a punctured 2-torus. 

In the case of $\Q_4$, it may be convenient to consider the cubical complex $\Q^\square_4$. 
The 6 vertices of  $\Q^\square_4$ correspond to the permutations $$1234, 2314, 3124, 2134, 1324, 3214.$$
Its 12 edges are 
 $$[12]34, [23]14, [13]24, 1[23]4, 3[12]4, 2[13]4, 12[34], 23[14], 31[24], 21[34], 13[24], 32[14].$$
and the three 2-cells are
$$[12][34], [23][14], [13][24].$$
Each 2-cell may be collapsed to a ``cross":
$$\includegraphics[width=200pt]{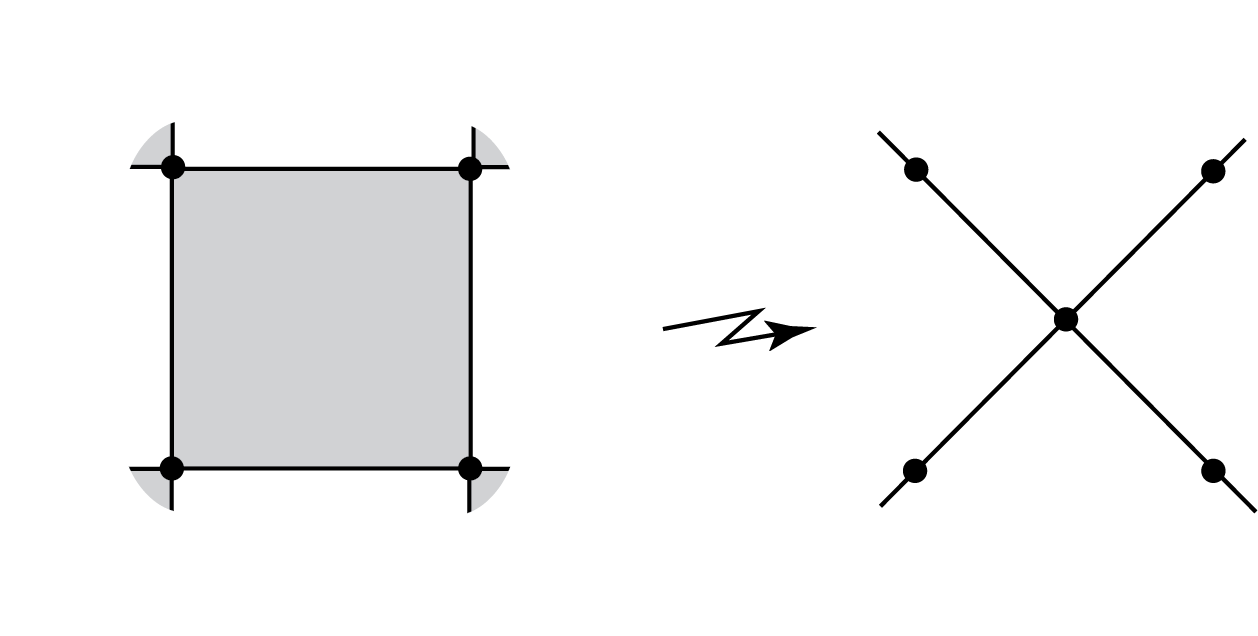}$$
After erasing the 6 bivalent vertices, we obtain a connected graph with 3 vertices and 6 edges; it is homotopy equivalent to a bouquet of 4 circles.

\end{proof}

\begin{theorem}
The space $\Q_5$ is homotopy equivalent to the Riemann surface of genus 4.
\end{theorem}
\begin{proof}
The cubical complex $\Q^\square_5$ is a connected 2-dimensional manifold. Indeed, every edge is adjacent to precisely two faces; for instance, the edge $[12]345$ is adjacent to $[12][34]5$ and $[12]3[45]$. Each vertex is adjacent to 5 edges and 5 faces; 
for instance, $12345$ is adjacent to the edges $[12]345$, $1[23]45$, $12[34]5$, $123[45]$ and $234[15]$ and the faces
$[12][34]5$, $[23]4[15]$,  $[12]3[45]$,  $2[34][15]$ and $[23]4[15]$.

According to Theorems~\ref{firstbetti} and \ref{euler}, the first Betti number of $\Q_5$ is 8 and the second Betti number is 1, so $\Q^\square_5$ is orientable of genus 4.
\end{proof}

Theorems~\ref{firstbetti}, \ref{euler} and \ref{cohdim} also give all the Betti numbers for $\Q_6$ and $\Q_7$:  
\begin{theorem}
We have $$b_1(\Q_6)=15,\quad b_2(\Q_6)=14,$$
$$b_1(\Q_7)=26,\quad b_2(\Q_7)=115,$$
and $b_i(\Q_6)=b_i(\Q_7)=0$ for $i>2.$
\end{theorem}

While $\Q_6$ and $\Q_7$ have cohomological dimension two, the cubical complexes $\Q^\square_6$ and $\Q^\square_7$ are 3-dimensional. In fact, one can construct an explicit 2-dimensional cubical complex homotopy equivalent to $\Q_6$:
\begin{theorem}\label{cubi}
There exists an $S_6$-equivariant 2-dimensional cubical complex with 30 vertices, 120 edges and 90 faces, homotopy equivalent to $\Q_6$. The action of $S_6$ is transitive on the vertices, edges and faces. 
\end{theorem}

\begin{proof}
The cubical complex  $\Q^\square_6$ has $120$ vertices, $360$ edges, 
$270$ 2-faces and $30$ 3-faces. 
\begin{itemize}
\item
Every vertex lies in the boundary of precisely two 3-faces. For instance, $123456$ lies in the boundary of $[12][34][56]$ and of $[23][45][16]$. 
\item
Each edge lies in the boundary of exactly three 2-faces. For example, $[12]3456$ lies in the boundary of  $[12][34]56$, $[12]3[45]6$ and $[12]34[56]$.  
\item
Out of the 270 2-faces, 180 lie in the boundary of one 3-face [the face $[12][34]56$ lies in the boundary of $[12][34][56]$] and 90 2-faces [such as  $[12]3[45]6$, for example)  in the boundary of no 3-face. 
\end{itemize}
Therefore, without changing the homotopy type of the complex, each 3-face can be collapsed on its subspace homeomorphic to the cone on the 1-skeleton of a cube: 
\smallskip

$$\includegraphics[width=250pt]{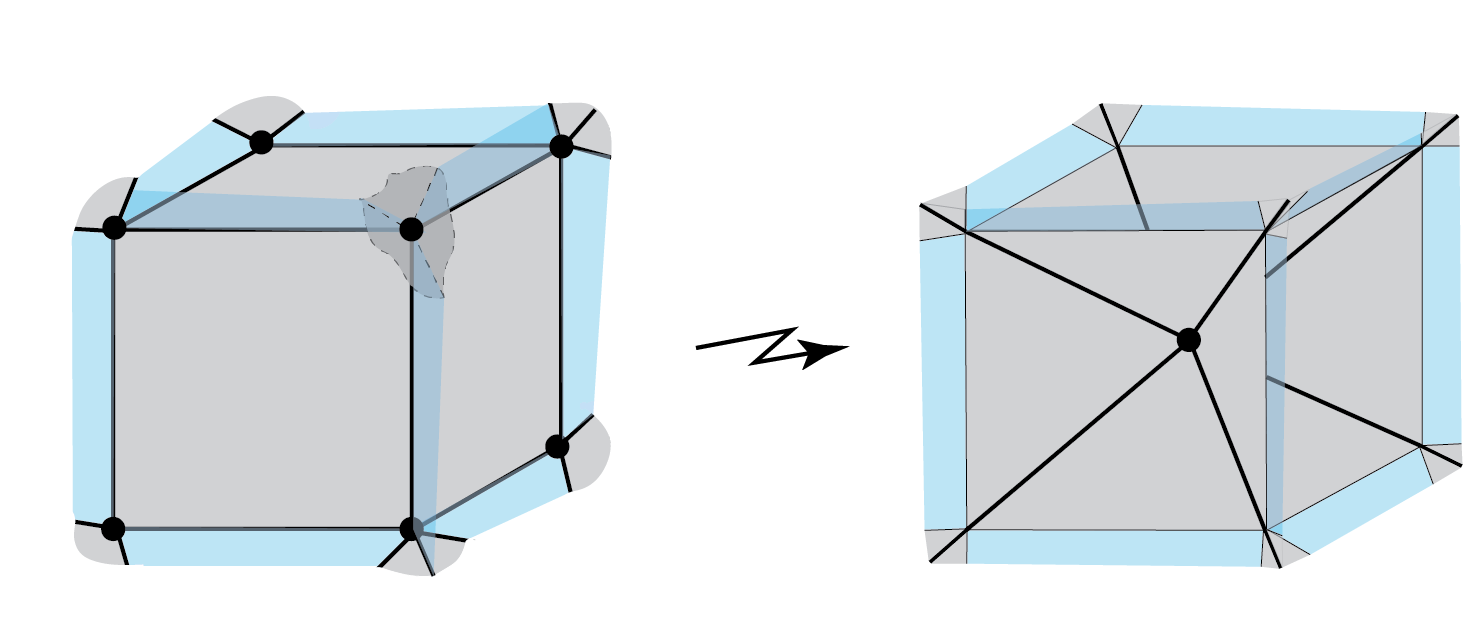}$$

\smallskip

The resulting 2-dimensional complex has 30 vertices, 120 edges and 90 square 2-faces.

There is a natural action of $S_6$ on $\Q^\square_6$ which is transitive on the vertices, edges, 3-faces and those 2-faces that do not lie in the boundary of any 3-face. The collapsing operation commutes with this action. The vertices, the edges and the 2-faces of the new complex are in 1-1 correspondence with the 3-faces, the vertices and the 2-faces not contained in the boundary of any 3-face of $\Q^\square_6$, respectively.
\end{proof}

\begin{remark}
In the 2-dimensional complex constructed in Theorem~\ref{cubi}, each vertex has a neighbourhood homeomorphic to the cone on the 1-skeleton of a cube. However, the stabilizer of a vertex in $S_6$ is not the full group of symmetries of a cube but is isomorphic to $(\Z/2)^3$. 
\end{remark}

\section{Round twin groups and pure cactus groups}

One motivation for the study of the round twin groups is their relation to the pure cactus groups, the fundamental groups of the moduli spaces  $\overline{\mathcal{M}}_{0,n}(\R)$ of stable real rational curves with marked points. We refer to \cite{Devadoss, Kapranov} for the definitions and properties of these spaces.

\begin{figure}[ht]
\includegraphics[width=1.8in]{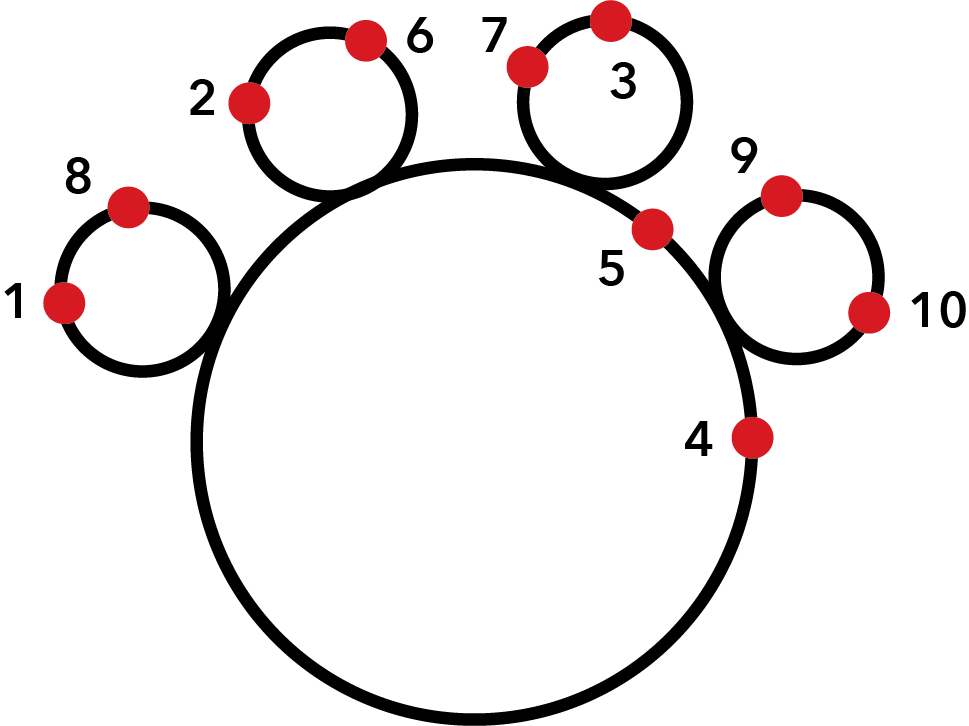}
\caption{A point in $\overline{\mathcal{M}}_{0,10}(\R)$ coming from $\Q_{10}$.}\label{cact}
\end{figure}

For $n\geq 3$, the moduli space ${\mathcal{M}}_{0,n}(\R)$ of real rational curves with $n$ marked points is the quotient of the configuration space of $n$ distinct points in $S^1$ under the action of the group $PSL(2,\R)$ of M\"obius transformations. This quotient map extends to a continuous map
$$q_n:\Q_n\to\overline{\mathcal{M}}_{0,n}(\R).$$
A configuration with points of multiplicity two is sent by $q_n$ to the class of a stable curve which has an extra component with three special points for each double point of the configuration; see Figure~\ref{cact}.

The map $q_n$, is actually a fibre bundle whose fibre may be identified with the stabilizer of $\infty$ in $PSL(2, \R)$.  This stabilizer has two contractible components. Therefore, we may think of $\Q_n$, up to homotopy, as a subset of the oriented double cover of $\overline{\mathcal{M}}_{0,n}(\R)$. 

\subsection{Moduli spaces of stable real rational curves with 5 and 6 marked points}

For $n=5$, every point of  $\overline{\mathcal{M}}_{0,5}(\R)$ is the image of a point in $\Q_5$, and  $\Q_5$ is homotopy equivalent to the double cover of $\overline{\mathcal{M}}_{0,5}(\R)$. This is consistent with our identification of 
$\Q_5$ as an orientable closed surface of genus 4 since $\overline{\mathcal{M}}_{0,5}(\R)$ is a connected sum of five projective planes.

\begin{figure}[ht]
\includegraphics[width=4in]{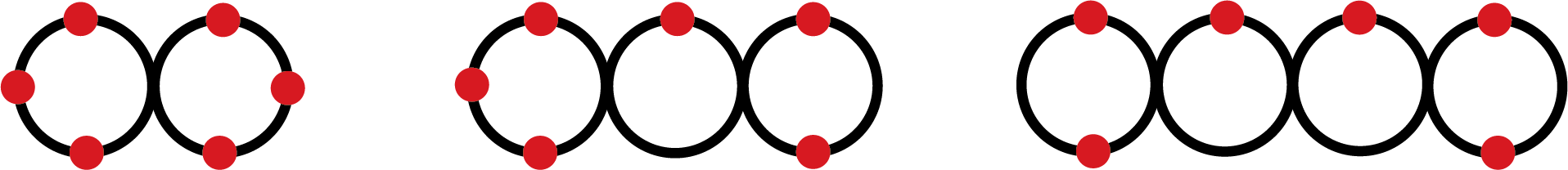}
\caption{Points in $\overline{\mathcal{M}}_{0,6}(\R)$ not coming from $M_{6}{(S^1)}$.}\label{badsix}
\end{figure}

When $n=6$, the points of $\overline{\mathcal{M}}_{0,6}(\R)$ that do not come from $\Q_6$ form the closure of the subset consisting of all curves with two components, each component carrying three marked points. The combinatorial types of all the curves in this closure are shown in Figure~\ref{badsix}.

In fact, $\overline{\mathcal{M}}_{0,6}(\R)$ is the blowup of $\R P^3$ along the configuration of points and lines that consists of the 4 vertices and the barycenter of a tetrahedron and all the lines passing through them.
 First, one blows up $\R P^3$ at the 5 points and then, at the 10 lines connecting them. 
The complement of $q_6(\Q_6)$ is then precisely the exceptional divisor of the blowup along the 10 lines, since this exceptional divisor consists of the curves on Figure~\ref{badsix}; see \cite{Devadoss}. 

This means that $q_6(\Q_6)$  is the complement to a 10-component link in the blowup of $\R P^3$ at 5 points.
As a consequence, the group $\X_6$ is the fundamental group of the complement of a 10-component link in the 3-manifold which is the boundary of a 4-disk with five 1-handles. It may be seen that its abelianization $H_1(\X_6)$ is generated by the classes of 1-handles and the meridians of the link components.



\subsection{The ``full'' round twin group} 
\begin{figure}[ht]
\includegraphics[width=3.5in]{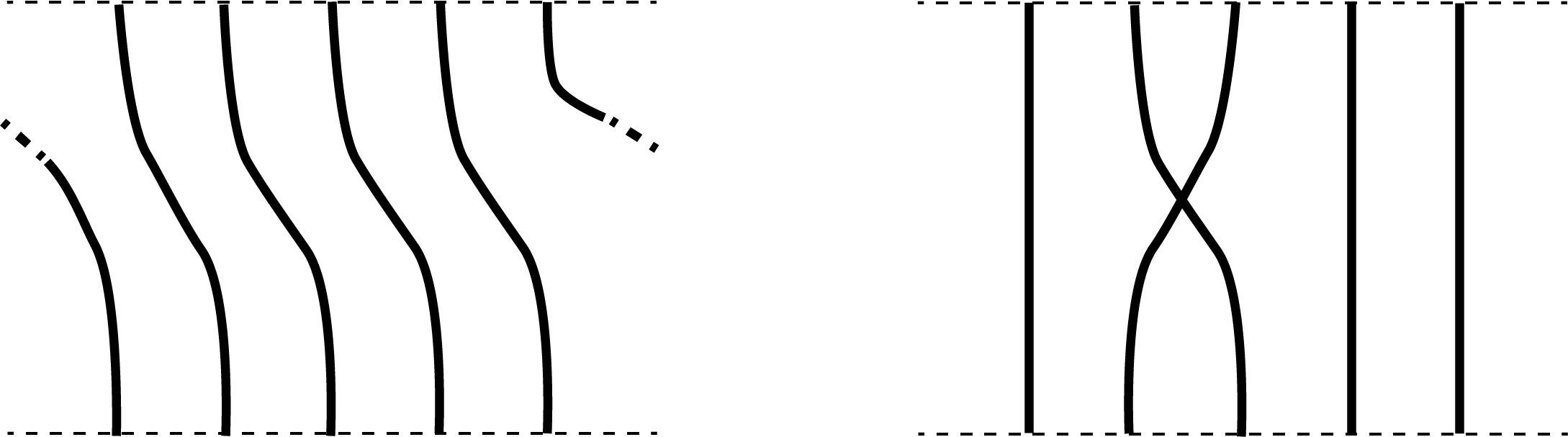}
\caption{The generators $\zeta$ and $\sigma_2$ of $\Y_5$.}\label{generatorsfull}
\end{figure}

Choose a basepoint $(x_1,\ldots, x_{n-1})\in \Q_n$ such that $x_i\neq x_j$  for all $i,j\in\{1,\ldots, n\}$ with $i\neq j$ (recall that $x_n=\infty$). Then, an element of $\X_n$ is determined by an $n-1$-tuple of paths in $S^1$, connecting each point $x_i$ ($i<n$) with itself. Together, these paths can be drawn as planar pure braids (twins) on $n-1$ strands that may ``pass through infinity'',  as in Figure~\ref{generatorsfull}. 

For a permutation $\tau$ of the set $1,\ldots, n-1$, consider an $n-1$-tuple of paths in $S^1$ connecting $x_i$ to $x_{\tau(i)}$ which, taken together with the constant path at infinity, determine a path in $\Q_n$. The fixed-end homotopy class of such a path in $\Q_n$ gives an element of what we call the  ``full'' round twin group $\Y_{n-1}$. We have a short exact sequence
$$1\to \X_{n}\to \Y_{n-1}\to S_{n-1}\to 1.$$

The group $\Y_{n-1}$ has a presentation with the generators $\sigma_1, \ldots, \sigma_{n-1}$ and $\zeta$ (see Figure~\ref{generatorsfull}) subject to the following relations:
\begin{equation}\label{relroundtwin}
\begin{array}{rcll}
\sigma_i^2&=& 1&\quad \text{for all\ } 1\leq i < n-1;\\
\sigma_i\sigma_j&=& \sigma_j\sigma_i &\quad \text{for all\ } 1\leq i, j < n-1 \text{\ with\  } |i- j| >1;\\
\sigma_i\zeta& = &\zeta\sigma_{i+1} &\quad \text{for all\ } 1\leq i < n-1.
\end{array}
\end{equation}

\subsection{The  round twin group as a subgroup of the pure cactus group} 

The group $$\Gamma_n=\pi_1 \overline{\mathcal{M}}_{0,n}(\R)$$ is known as the $n$th \emph{pure cactus group}.
The map $q_n:\Q_n\to\overline{\mathcal{M}}_{0,n}(\R)$ gives rise to a homomorphism
$$\Pi_n\to \Gamma_n$$
for each $n$. While we do not have neat presentations for either of the two series of groups, this homomorphism can be described very explicitly in terms of the corresponding full groups. 

The \emph{full} cactus group $J_n$ (see \cite{HK}) has a presentation with the generators $s_{p,q}$, where $1\leq p< q\leq n$, and the following relations:
\begin{equation}\label{eqcac}
\begin{array}{rcll}
s_{p,q}^2&=& 1,&\\
s_{p,q}s_{m,r}&=& s_{m,r}s_{p,q} &\quad \text{if\ } [p,q]\cap [m,r] =\emptyset,\\
s_{p,q}s_{m,r}&=& s_{p+q-r, p+q-m}s_{p,q} &\quad \text{if\ }  [m,r] \subset [p,q].
\end{array}
\end{equation}
There is a homomorphism $J_n\to S_n$ to the symmetric group: it sends $s_{p,q}$ into the permutation $\tau_{p,q}$ of $\{1, \ldots, n\}$ which reverses the order of $p, p+1,\ldots, q$ and leaves the rest of the elements unchanged. The {pure cactus group} $\Gamma_{n+1}$ is the kernel of this homomorphism.

Consider the homomorphism
$$\kappa: \Y_n\to J_n$$
defined by
$$\begin{array}{lcl}
\kappa(\sigma_i) &= &s_{i,i+1},\\
\kappa(\zeta) &= &s_{1,n} s_{2,n},
\end{array}
$$
with $\sigma_i$ and $\zeta$ as in (\ref{relroundtwin}).
This homomorphism  sends $\Pi_{n+1}\subset  \Y_n$ to $\Gamma_{n+1}$. We will denote by the same letter the map sending words in the generators of  $\Y_n$ to words in the generators of $J_n$.

In order to see that the restriction of $\kappa$ to  $\Pi_n$ is actually induced by the map $$\Q_n\to \overline{\mathcal{M}}_{0,n}(\R),$$ one has to recall the geometric meaning of the generators of $J_n$:
the generator $s_{p,q}$ corresponds to a path in  $\overline{\mathcal{M}}_{0,n+1}(\R)$ in which the marked points $p,\ldots, q$ collide and bubble off onto a new component
and then return to the original component in the reversed order (see \cite{HK}). This shows that the $\sigma_{i}$ map into the
$s_{i,i+1}$ and $\zeta$ must go to $s_{1,n} s_{2,n}$.

\medskip

\begin{theorem}\label{inj}
The homomorphism $\kappa:\Y_n\to J_n$ is injective for $n\geq 4$.
\end{theorem}
In particular, for $n\geq 4$, the round twin groups embed into the corresponding pure cactus groups. Theorem~\ref{inj} fails for $n<4$:  for instance, $\Pi_4=F_4$ while $\Gamma_4$ is infinite cyclic.

This statement has the following corollary:
\begin{theorem}
The group $\X_n$ is residually nilpotent for all $n$.
\end{theorem}
Indeed, from \cite{MCact} we know that the pure cactus groups are all residually nilpotent.

\begin{proof}
If $w$ is a word in the generators $s_{p,q}$ that defines the trivial element of $J_n$, there exists a sequence $w_1,\ldots, w_m$ such that $w_1=w$, $w_n$ is trivial, $w_{i+1}$ is obtained from $w_i$ by applying one of the relations (\ref{eqcac})
once, and the length of $w_{i+1}$ is not greater than the length of $w_i$.

This follows from the proof of Proposition~2 in \cite{MCact} where it is shown that two words in the $s_{p,q}$ which represent the same element of $J_n$ and are \emph{locally reduced} (their lengths cannot be decreased by applying the relations (\ref{eqcac})) must have the same length. Indeed, if $w$ can only be taken into the trivial word by a sequence of moves that increases the length at some point, there exists a locally reduced word which represents the trivial element in $J_n$, which is impossible.

The image of a word in the generators $\sigma_i$ and $\zeta$ under $\kappa$ is a word $w$ in the $s_{i,i+1}$, $s_{1,n}$ and $s_{2,n}$.  If it represents the trivial element of $J_n$, it can be transformed into the trivial word by means of the relations that involve only the generators $s_{i,i+1}$, $s_{1,n}$, $s_{2,n}$ and $s_{1,n-1}$, since the appearance of any other generator of $J_n$ in the sequence of the words connecting $w$ and $1$ would imply that at some point the length of the word increases. 

Assume that $n\geq 4$; under this condition neither of $s_{1,n}$, $s_{1,n-1}$ or $s_{2,n}$ coincides with any of $s_{i,i+1}$. Let $z$ be the word $s_{1,n} s_{2,n}$. For any word $u$ in the generators $s_{i,i+1}$, $s_{1,n}$, $s_{2,n}$ and $s_{1,n-1}$, define the word $\mu(u)$ in the $s_{1,2}$,  $\ldots,$ $s_{n-1,n}$, $z$ and $s_{1,n}$ inductively as follows. 

For a word $w$ in the generators $s_{p,q}$, let $\overline{w}$ be the word in which each $s_{p,q}$ is replaced by $s_{n-q+1,n-p+1}$. Now: 
\begin{itemize}
\item if $u=1$ we set $\mu(u)=1$;
\item if $u=s_{i,i+1} v$, we set $\mu(u)=s_{i,i+1} \mu(v)$;
\item if $u=s_{1,n} v$, we set $\mu(u)=\overline{\mu(v)} s_{1,n}$;
\item if $u=s_{2,n} v$, we set $\mu(u)=z^{-1} \overline{\mu(v)} s_{1,n}$;
\item if $u=s_{1,n-1} v$, we set $\mu(u)=z \overline{\mu(v)} s_{1,n}$. 
\end{itemize}

For any word $u$ in  $s_{i,i+1}$, $s_{1,n}$, $s_{2,n}$ and $s_{1,n-1}$, the word $\mu(u)$ is of the form $\mu'(u)  s_{1,n}^k$, where $\mu'(u)$ is a word in $s_{i,i+1}$, and $z$ only.

Assume $v$ is a word in the generators  $\sigma_i$ and $\zeta$ such that $\kappa(v)$ defines a trivial element of $J_n$.
Take the sequence of words  $w_1,\ldots, w_n$ in the $s_{p,q}$, such that $w_1=\kappa(v)$, $w_m$ is trivial, $w_{i+1}$ is obtained from $w_i$ by applying one of the relations (\ref{eqcac}) once and the length of $w_{i+1}$ is not greater than the length of $w_i$. Then the sequence of words  $\mu'(w_1),\ldots, \mu'(w_n)$ in $s_{1,2}$,  $\ldots,$ $s_{n-1,n}$ and $z$ transforms $\kappa(v)$ into the trivial word by means of the relations (\ref{eqcac}). Replacing each $s_{i,i+1}$ with $\sigma_{i}$ and $z$ with $\zeta$, we obtain a sequence of words which transforms $v$ into the trivial word by means of the relations (\ref{relroundtwin}) in $\Y_n$. In particular, this means that $v=1$ in $\Y_n$. 
\end{proof}

\begin{remark}
The fact that the twin groups inject into the cactus groups has been observed in \cite{BCL}.
\end{remark}

\end{document}